\magnification=\magstep1

\def\Z{\bf Z}
\def\g{\Gamma}
\def\ga{\Gamma^{op}}
\def\s{\Sigma}
\def\N{{\bf N}}
\def\gs{{\cal G}{\cal S}}
\def\gss{{\cal G}{\cal S}{\cal S}}
\vskip1cm

\centerline{\bf  MULTIPLICATIVE MAPS FROM $H\Z$  TO A RING
SPECTRUM $R$}

\centerline{\bf - A NAIVE VERSION}

\medskip

\centerline{\it by Stanislaw Betley}

\vskip1cm

 {\bf 0. Introduction.}

\bigskip
For a commutative ring $B$ Stefan Schwede described in [S] the
surprising connection between stable homotopy theory of
commutative $B$-algebras and formal group laws over $B$. The
stable homotopy operations of commutative simplicial B-algebras
are described by the algebra $\pi_* DB$, where $DB$ is a certain
''classical'' spectrum studied by Bousfield, Dwyer and others, see
for example [D]. Schwede was able to describe the weak homotopy
type of the space of multiplicative  maps from the
Eilenberg-MacLane spectrum $H\Z$ to $DB$ in terms of formal group
 laws over $B$ and their isomorphisms. But his methods seem to be much more general and
 should work in other situations as well. On the other hand the formal
 group laws over $B$ are present in the description of $DB$ so, in order to
 generalize his results one should start from defining something like
 ''formal group law'' even with the lack of formal power series.

 In the present note we offer a definition of a formal group law
 in a ring spectrum $R$. With it we recover the weak version of the $\pi_0$-result of
 Schwede with any ring spectrum $R$ at the place of $DB$. The
 obvious generalization of the full Schwede's result is clearly
 visible but we don't have any evidence to call it even
 ''conjecture''. At present we do not see methods of attacking
 this problem in full generality.

 Every multiplicative map $H{\Z} \to R$ give $R$ the structure of a
 ring over $H\Z$. We hope that  observations presented in this note can be fruitful
 for our better understanding of the category of such objects.

 We use here the language of Lydakis from [L], which we summarize in Section 1. The ring spectrum
 means here a $\g$-ring in the sense of  [L].
 It means that a ring spectrum is a functor from the category of
 finite sets to simplicial sets with the extra structure. Maps
 between such objects are given in terms of natural
 transformations of functors. The word ''naive'' in our title
 refers to the fact that we work mostly with combinatorial structures only, so we don't have
 to use (up the last two  pages) models for our spectra
 which give the correct homotopy type of the mapping space
 (fibrant - cofibrant replacement). We would like to thank the
 referee for many useful suggestions which helped us to improve
 the presentation. Especially the names used in the definitions
$ 2.1$ and $3.1$ should be blamed on her/him.

 \vskip2cm

{\bf Acknowledgment:} This research was partially supported by the
Polish Scientific Grant N N201 387034.

\vfill
\eject

{\bf 1. Preliminaries on $\g$-spaces and $ \g $-rings.}

\bigskip

Let  $[n]$ denote the pointed set $\{ 0, 1,...,n\}$ with $0$ as a
basepoint for a nonnegative integer $n$.

We want to distinguish three types of pointed set maps which will
play the crucial role in the future. Two of them map $[n]\to
[n-1]$ and the third one goes the other way around. The map
$p_i^n:[n]\to [n-1]$ defined by $p_i^n(j)=j$ for $j<i$,
$p_i^n(i)=0$, $p_i^n(j)=j-1$ for $j>i$ will be called the ith
restriction. For any $i<j\leq n$ we have the ''summing'' map
$s_{i,j,k}^n:[n]\to [n-1]$ defined via the formula
$s_{i,j,k}^n(i)=s_{i,j,k}^n(j)=k$, and the other elements $a\in
[n]$ are mapped bijectively onto $[n-1]\setminus {k}$, preserving
ordering. The third map $d^n_j:[n-1]\to [n]$ is injective, misses
$j\in [n]$ and preserves the order.

The category $\ga$ is a full subcategory of the category of
pointed sets, with objects all $[n]$.  The category of $\g$-spaces
is the full subcategory of the category of functors from $\ga$ to
pointed simplicial sets with objects satisfying $F[0]=[0]$ and
maps given by the natural transformations of functors. Perhaps one
should explain here that the notation $\ga$ comes from the fact
that our category is dual to Segal's category $\g$ from [Se].
Every $\g$-space can be prolonged by direct limits to the functor
defined on the category of pointed sets. In our notation we will
not distinguish between the $\g$-space and the described above
extension. We will use capital letters $K$, $L$, ... for denoting
pointed sets. In the future, if we need ordering of the pointed
set  $[n]\wedge [m]$ which identifies it with  $[nm]$  we will use
always the inverse lexicographical order.

\medskip

{\bf Convention:} If it causes no misunderstandings having pointed
sets $K$ and $L$ and a  pointed map $f:K\to L$ we  write $f$
instead of $F(f)$ for the induced map $F(K)\to F(L)$.

\medskip

For a $\g$-space $F$ let $RF$ denote the $\g \times \g$ space
defined as $RF(K,L)=F(K\wedge L)$.  Having two $\g$-spaces $F$ and
$F'$ we can form their exterior smash product $\g \times \g$-space
$F\tilde\wedge F'$ which is defined by the formula $F\tilde\wedge
F'(K, L)=F(K)\wedge F'(L)$. Then the smash product of $F$ and $F'$
is the universal $\g$-space $F''$ with a map of $\g \times
\g$-spaces $F\tilde\wedge F'\to RF''$, (see [L, Remark 2.4]).
Moreover, if we denote by $\gs$ the category of $\g$-spaces and by
$\gss$ the category of $\g \times \g$-spaces then for given
$\g$-spaces $F_1$, $F_2$ and $F_3$ we have (following [L,Theorem
2.2]):
$$\gs(F_1\wedge F_2,F_3)=\gss(F_1\tilde\wedge F_2,RF_3)$$

\medskip

{\bf Remark 1.1:} The symmetric group $\s_n$ acts on the set $\{
0,1,...,n\}$ by permuting $\{ 1,...,n\}$ and hence this group acts
on $F[n]$ for any  $\g$-space $F$. We will use this action
restricted to various subgroups of $\s_n$ in the future.

\medskip

Let ${\bf S}$ denote the $\g$-space defined by the identity
functor. We say that $\g$-space $F$ is a  $\g$-ring if there are
maps $\eta :{\bf S}\to F$ called the unit and $\mu: F\wedge F\to
F$ called the multiplication satisfying usual associativity
 and unit conditions (see [L, 2.13]).

\bigskip

{\bf Remark 1.2:} By our previous
 observations  $\mu$ is determined by a map
 $\tilde\mu:F\tilde\wedge F \to RF$, which is fully determined  by a collection of maps  $\tilde\mu:F[n]\wedge F[m] \to
 F([n]\wedge [m])$ natural in $[n]$ and $[m]$  and
 satisfying obvious associativity conditions.

\bigskip

Let us introduce one more piece of notation. We will say that a
$\g$-ring $R$ is discrete if for any pointed set $K$ , $R(K)$ is
just a set considered as a simplicial set in a trivial way. Assume
that $R$ is a discrete $\g$-ring. Then $R[1]$ is a unital monoid
with zero. Moreover $\eta$ takes $1\in {\bf S}[1]$ to the unit of
$R[1]$.

\bigskip

{\bf Remark 1.3:} Assume that $R$ is a discrete $\g$-ring. Then
the map $\tilde\mu:R(K)\wedge R(L)\to R(K\wedge L)$ of 1.2 is a
map of sets which is associative with respect to the smash product
of pointed sets. It means that if $p\in R(K)$ and $q\in R(L)$ then
it makes sense to say that the product of $p$ and $q$ belongs to
$R(K\wedge L)$ which, of course , means that $\tilde\mu (p,q)\in
R(K\wedge L)$. We will usually write the result of such
multiplication as $pq\in R(K\wedge L)$.

\bigskip

{\bf 2. Multiplicative maps from $H\N$ to a discrete $\g$-ring
$R$}

\bigskip

We want to study multiplicative maps from the Eilenberg-MacLane
spectrum $H\Z$ to a $\g$-ring $R$. But the answer is a bit
technical and we will postpone it until section 3. In the present
section we will study maps from the spectrum stably equivalent to
$H\Z$ which is easier to study. Let us start from recalling a
$\g$-ring model of $H\Z$. As a functor $H\Z$ takes $K$ to a
reduced free abelian group generated by $K$. The map

$$\eta:{\bf S} \to H{\Z}$$

\noindent is given by the embedding of generators. The
multiplication map
$$\mu: \tilde{\Z}(K)\wedge \tilde{\Z} (L) \to \tilde{\Z}(K\wedge L)$$
is defined via the formula

$$(\sum_{k\in K}a_kk )\wedge (\sum_{l\in L}b_ll) \mapsto \sum_{k\wedge l\in
K\wedge L} a_kb_l(k\wedge l)$$

The $\g$-ring  $H\N$ is defined with the same formulas but with
the additive monoid of natural numbers (with $0$ ) instead of
integers. The embedding $H\N\to H\Z$ induces a stable equivalence
of spectra because for any $k>0$ the map $H\N(S^k)\to H{\Z}(S^k)$
is a homotopy equivalence by a theorem of Spanier [Sp, Theorem
4.4]. This map is obviously multiplicative but there is no
nontrivial multiplicative map going the other way.

It is easy to realize that Schwede's map $H{\Z} \to DB$ associated
to a formal group $F$ is determined uniquely by saying that the
image of $(1,1)\in H{\Z}[2]$ is $F$. We comment on this more later
but now we should  define the possible images of $(1,1)\in H\N[2]$
in the case of an arbitrary $\g$-ring $R$. Below we write the
first definition of a formal group law in a $\g$-ring $R$.

\bigskip

{\bf Definition 2.1:} A  formal sum law in a $\g$-ring $R$ is an
element $w\in R[2]$ which  satisfies the following properties:

1. $p^2_1(w)=1$, $p^2_2(w)=1$,

2. any power $w^{k}\in R[2^k]$ is fixed under the action of the
symmetric group $\Sigma_{2^{k}}$.

\bigskip

{\bf Theorem 2.2:} Let $R$ be a discrete $\g$-ring. Then every
formal sum law in $R$ determines the multiplicative map $\phi: H\N
\to R$.

\medskip

Proof. Let $w$ be a  formal sum law in $R$. Let $1_n=(1,...,1)\in
H\N[n]$. We want to show that associating  $ 1_{2^n}\mapsto w^n$
defines the desired map $\phi$.

Observe first that any element $(n_1,...,n_k)\in {\N} [k]$ can be
presented as an image of $ 1_n$ for a certain $n$ . So our map is
uniquely determined on elements $ 1_n$: if for a pointed map
$f:[n]\to [k]$ we have  $f(1_n)=(n_1,...,n_k)$ then we must have

$$\phi(n_1,...,n_k)=f(w^n)$$

Hence  we only have to show that $\phi$ is well defined by the
formula above.

First of all, for the given $k$tuple $(n_1,...,n_k)\in {\N} [k]$
there is a minimal $n$ such that $ 1_{2^n}$ maps to
$(n_1,...,n_k)$ . Obviously $n$ is equal to the minimal natural
number $n$ which satisfies the  condition $2^n\geq
\sum_{i=1}^kn_i$. Of course there are many ways of mapping $
1_{2^n}$ to $(n_1,...,n_k)$ but all of them give the same
definition of  $\phi (n_1,...,n_k)$
 because of the condition $2$ of the definition of the formal sum
 law.

 Assume now that $g( 1_{2^m})=(n_1,...,n_k)$ for a certain map $g$ with
   $n<m$. Then it is easy to see that $g$ factors through $ 1_{2^n}$.
  Hence our proof that $\phi$ is well defined will be finished if
  we show:
\medskip

{\bf Lemma 2.2.1:} For any $k$,  $w^{2^k}$ maps to $w^{2^{k-1}}$
under any map $f$ which satisfies $f(1_{2^k})=1_{2^{k-1}}$.

\medskip

Proof of 2.2.1. Assume first that  $f$ takes last $2^{k-1}$
coordinates to zero. In other words, using the fact that
$$[2^k]=[2^{k-1}]\wedge [2]$$

we can write

$$f= Id_{[2^{k-1}]}\wedge  p_2^2$$

Then

$$f(w^{2^k})=
f(w^{2^{k-1}}\cdot w)=(Id_{[2^{k-1}]}\wedge
p_2^2)(w^{2^{k-1}}\cdot w)=w^{2^{k-1}}\cdot 1= w^{2^{k-1}} $$

Now observe that, because of the property $2$ of  formal sum laws,
it is enough to consider only maps like $f$ above. Any other $f'$
which takes $ 1_{2^k}$ to $ 1_{2^{k-1}}$ differs from $f$ by an
action of an element from $\Sigma_{2^k}$.

We will finish the proof of the theorem if we show that our map
$\phi$ obtained in the way described above is multiplicative. To
check this observe first that if $f(1_{2^n})=(n_1,...,n_k)$ and
$g(1_{2^m})=(m_1,...,m_l)$ then $(f\wedge g)(1_{2^n}\cdot
1_{2^m})=(f\wedge g)(1_{2^{n+m}})=(n_1,...,n_k)(m_1,...,m_l)$ as
elements of $H\N [kl]$. We calculate further:

$$\phi((n_1,...,n_k)\cdot (m_1,...,m_l)) =\phi(f\wedge g(1_{2^{n+m}}))=f\wedge g\phi(1_{2^{n+m}})=f\wedge
g(w^n\cdot w^m)=$$

$$=f(w^n)\cdot g(w^m)=\phi(f(1_{2^n}))\cdot \phi (g(1_{2^m}))=
\phi(n_1,...,n_k)\cdot \phi(m_1,...,m_l)$$

and the proof is finished.

\bigskip

We will come to the issue when two formal sum laws give homotopic
maps later in a more general setting. But from the proof of
theorem 2.2 it is easy to derive the following observation:

\bigskip

{\bf Remark 2.3:} Observe that our assumption that $R$ is discrete
is not important. We could define a  formal sum law in $R$ as a
$0$-simplex  of $R[2]$ and the rest would go through by the same
arguments.

\bigskip

{\bf 3. Multiplicative maps from $H\Z$ to a discrete $\g$-ring
$R$.}

\bigskip

Now we move towards studying  multiplicative maps $H{\Z}\to R$. We
would like to define the formal group laws in this situation in
such a way that we get the same statement as in 2.2. But first of
all let us identify the complications which occur when we allow
negative coordinates in our source $\g$-ring. The problem is that
while working with  an arbitrary $\g$-ring $R$ one does not have
any natural way of defining maps coming from multiplying one
''variable'' by $-1$. That was not a problem in the case of $DB$.
More generally  this is not a problem in the case of any $\g$-ring
coming from the composition of functors
$$T\circ L:{\g}^{op} \to Sets_*$$

where $L$ is the  linearization functor from sets to the category
$B_{free}$ of  free modules over some ring $B$ and $T:B_{free}\to
Sets_*$. We plan to study such situations in the forthcoming paper
 but now we would like to define the formal group law in full
generality overcoming the difficulty described above.

 But before the definition
we have to describe a particular type of an action of
$\Sigma_{2^{k-1}}\times \Sigma_{2^{k-1}}$ on $F[2^k]$ for any
$\g$-space $F$. This action will ba called {\it special} later on.
Let for any $k$, $\pm 1_k$ be equal to $(1,-1)^k\in H{\Z}[2^k]$.
Our convention on ordering smash products of pointed sets permits
to split $[2^k]=A_+\vee A_{-}$ accordingly to the rule that
 $(1,-1)^k$ has $1$ at the coordinates from $A_+$ and $-1$ at
 $A_-$. In a more formal way we can say that an element $i\in \{
 1,...,2^k\} $ belongs to $A_+$ it and only if the binary
 expansion of $i-1$ has an even number of digits ''1''. There is
 also a ''coordinate-free'' way of describing the splitting. If we
 identify the set $[2^k]$ with the set of subsets of the set $\{ 1,2,...,k\}$ then
 $A_+$ ($A_-$) consists of sets of even (odd) order.

The special action of $\Sigma_{2^{k-1}}\times \Sigma_{2^{k-1}}$ on
$F[2^k]$ is defined as follows: if $a\times b \in
\Sigma_{2^{k-1}}\times \Sigma_{2^{k-1}}$ then $a$ permutes the
coordinates from $A_+$ and  $b$ permutes the rest of the
coordinates. Let $\sigma$ be the nontrivial element of $\Sigma_2$.

\bigskip

{\bf Definition 3.1:} A formal difference law in a $\g$-ring $R$
is an
 element $r\in R[2]$  satisfying the following
properties:

1. $p^2_2(r)=1$, $s^2_{1,2,1}(r)=0$

2. $p^2_1(r)r=rp^2_1(r) =\sigma (r)$ in $R[2]$

3. any power $r^{k}\in R[2^k]$ is fixed under the special action
of $\Sigma_{2^{k-1}}\times \Sigma_{2^{k-1}}$.

4. for any $k$ , $i<j$ and $l$  such that $i\in A_+$ and $j\in
A_-$ or $j\in A_+$ and $i\in A_-$ we have
$$s^{2^k}_{i,j,l}(r^k)= d^{{2^k}-1}_lp^{{2^k}-1}_ip^{2^k}_{j}(r^k)$$

\bigskip

 Observe first that $p_1^2(r)$ plays a role of $-1$ in $R[1]$
because we can calculate;

$$(p_1^2(r))^2=(p_1^2\wedge p_1^2)(r^2)=(p^2_1\wedge p^2_1)\circ
\tau )(r^2)=(p_2^2\wedge p_2^2)(r^2)=(p_2^2(r))^2=1$$

\noindent where $\tau$ is a special permutation in $\Sigma_2\times
\Sigma_2$ given by the transposition $(1,4)$. We can now  compare
our new definition  with results and definitions from Section 2.
We check that if $r$ is a formal difference law in a $\g$-ring $R$
then $w=p_2^3\circ p_2^4(r^2)\in R[2]$ is a  formal sum law in the
sens of definition 2.1. Indeed:

$$p^2_2(w)= p^2_2\wedge
p_2^2(r^2)=(p_2^2(r))^2=1$$

\noindent Similarly:
$$p^2_1(w)= p^2_1\wedge
p_1^2(r^2)=(p_1^2(r))^2=1$$

\noindent Moreover

$$\sigma (w))=p_2^3\circ p_2^4(\tau (r^2))=p_2^3\circ
p_2^4(r^2)$$

\noindent and hence $w$ is fixed under the action of $\Sigma_2$.
Let $p$ denote $p^3_2\circ p^4_2$.  By naturality of the smash
product and multiplication maps we have the following commutative
diagram:

$$\matrix{ & R[4]  \wedge  R[4] & \longrightarrow & R[16] \cr
            &\downarrow  &\   &\downarrow \cr
           &R[2] \wedge R[2] &\longrightarrow & R[4]  \cr}
  $$

\noindent where the  left vertical arrow is given by $p\wedge p$
and horizontally we have multiplication maps. Then the right
vertical map is defined by the set map which takes $4$ elements of
$A_{+}$ bijectively to nonzero elements of $[4]$ and the rest
elements of $[16]$ to $0$ . Hence the action of any permutation
from $\Sigma_4$ on $w\wedge w$ rises to the special permutation
acting on $R[16]$. This argument generalizes  easily to higher
degrees because $p^{\wedge k}$ maps bijectively $2^k$ elements of
 $A_{+}\subset [4^k]$ to nonzero elements of $[2^k]$ and has value
 $0$ otherwise.

Thinking about our definition as if it was a definition of a
formal group law in an ordinary sense we can give an
interpretation of the most of the structure described in 3.1. An
element $p_1^2(r)$ plays a role of $-1$ in the ''commutative group
structure'' defined by $r$. Hence it commutes ''with other
elements''. Condition 1 is always included in the general
definition of a formal group law. The same can be said about
condition 3 - in the classical case of formal power series this
kind of invariance property is indirectly in the definition of a
formal group law.

The condition 4 is new and makes the situation technically more
complicated.  It is hard to imagine what would be the abstract
meaning of it. This condition  is strongly related to the fact
that $1+(-1)=0$ in $\Z$ which is a very additive condition, having
no meaning in the structure of an arbitrary $\g$-ring . The
simplest explanation which one can imagine for the need of 4 is
the following remark:
 condition 4  is an extension of the second formula from condition
1 to higher degrees, which is needed when we face the lack of
additivity. The good news is that condition 4 from 3.1 is  often
satisfied in interesting cases, namely in $\g$-rings coming from
algebraic theories. This is the case of the $\g$-ring $DB$. We are
not going to define here what is an algebraic theory and what is
the definition of a  $\g$-ring associated to it. Instead we send
the interested reader to [S2, Section 2].

\bigskip

{\bf Remark 3.2:} Let $T^s$ be a $\g$-ring associated to the
algebraic theory $T$. Then condition 4 of 3.1 is satisfied for
$T^s$ as a consequence of condition 1.

\medskip

Proof (the sketch): We will follow [S1, Section 2] without further
explanations. Observe first that formula 4 of 3.1 for $k=1$ is
equivalent to the second equality of the condition 1 of 3.1 and
hence is satisfied. By definition $T^s[n]=hom_T([n],[1])$ and the
multiplication
$$T^s[n]\wedge T^s[m]\to T^s[nm]$$

\noindent is obtained from composition. It means that we can write
it as

$$\alpha\wedge \beta \mapsto \beta\circ (\alpha,...,\alpha)$$
with our convention of identifying $[n]\wedge [m]$ with $[nm]$. So
in the notation as above  $r^k=r^{k-1}\circ (r,...,r)$ and the
value of the  map $s^{2^k}_{i,j,l}$ on $r^k$ is the same as if we
apply $s^{2}_{1,2,1}$ to one of the $r$'s in the bracket by
naturality and property 3 of 3.1. So

$$s^{2^k}_{i,j,l}(r^k)=r^{k-1}\circ (r,...,r,0,r,...,r)= d^{{2^k}-1}_lp^{{2^k}-1}_ip^{2^k}_{j}(r^k)$$

\bigskip

{\bf Definition 3.3:}  A {\it homomorphism} $a:r_1\to r_2$ of
formal difference  laws in $R$ is an element $a\in R[1]$
satisfying $ar_1=r_2a$. An invertible homomorphism is called an
{\it isomorphism} of formal difference  laws.   An isomorphism is
called {\it strict} if it maps  to the unit component of $R$ under
the map $R[1]\to \pi_0R$.

\bigskip

Perhaps for completeness it is worth here to recall the definition
of the map $R[1]\to \pi_0R$. According to [S2, lemma 1.2] $\pi_0R$
can be presented as the cokernel of the map
$$ \tilde Zp^2_2 + \tilde Zp^2_1 - \tilde Zs^2_{1,2,1}:\tilde
Z[R[2]]\to \tilde Z[R[1]]$$

\noindent Then our map can be described as an embedding of
generators composed with the quotient map described above.

\bigskip

{\bf Theorem 3.4:} Let $R$ be a discrete $\g$-ring. Then every
formal difference  law in $R$ determines the multiplicative map
$\phi: H{\Z} \to R$.

\medskip

Proof. Let $r$ be a formal difference law in $R$. We want to show
that associating  $\pm 1_{n}\mapsto r^n$ defines the desired map.

Observe first that any element $(n_1,...,n_k)\in {\Z} [k]$ can be
presented as an image of $\pm 1_n$ for certain $n$. Hence our map
is uniquely determined on elements $\pm 1_n$ and we only have to
show that $\phi$ is well defined.

We would like to follow the proof of 2.2 but the situation is
different now. The proof of 2.2 was based on the fact that an
element $(n_1,...,n_k)\in {\N} [k]$ was equal to the image of
$1_n$ in a unique way up to a permutation  which acted  trivially
on the corresponding power of $r$. This is now not the case: $1$
and $-1$ from different coordinates in $\Z [k]$ can cancel either
by mapping coordinates to the base point or by the summing map. In
the case of the proof of 2.2 we had only to consider the first
possibility.

 First of all, as previously, for the given $k$tuple
$(n_1,...,n_k)\in {\Z} [k]$ there is minimal $n$ such that $\pm
1_n$ maps to $(n_1,...,n_k)$ by the map $f'$. We can assume that
all $n_i$s are different from $0$.  There is a special permutation
$\sigma$ such that $f=f'\circ \sigma$ takes first $\vert
n_1\vert$-coordinates in $\sigma ((\pm 1)^k)$ with the same sign
as $n_1$ to $n_1$, next $\vert n_2\vert$ coordinates with correct
signs  to $n_2$ and so on. Let $N_k=\Sigma_{i=1}^k \vert
n_i\vert$. It means  that all ones and minus ones on the other
$n-N_k$ coordinates  have to cancel to zero. Assume that $a<b$ and
we have 1 on $a$th coordinate and $-1$ on $b$th and they cancel
each other (add to 0). Then obviously $f=f\circ d^{2^n}_b\circ
s_{a,b,a}^{2^n}$ as the maps of pointed sets and we can iterate
this process composing with more pairs of maps $d^{2^n}_*\circ
s_{*,*,*}^{2^n}$. But observe that, because of condition 4 of the
definition of the formal difference  law we have

$$f(r^n)=f\circ d^{2^n}_b\circ s_{a,b,a}^{2^n}(r^n)=f\circ
d^{2^n}_b\circ d^{2^{n}-1}_a\circ p^{{2^n}-1}_a\circ
p^{2^n}_{b}(r^n)$$

It means that the value of $f$ on $r^n$ is the same as the value
of a map which takes $a$ and $b$ to a base point. Iterating this
process we see that we have justified the lemma:

\medskip

{\bf Lemma: 3.4.1:} Let $g:[2^n]\to [k]$ be a map which agrees
with $f$ on the  $N_k$ elements chosen as it is described above
and takes the rest to the base point. Then

$$f(r^n)=g(r^n)$$

\medskip
Hence we see that our map $\phi$ is well defined. Whichever map
$f$ taking $\pm 1_n$ to $(n_1,...,n_k)$ we use it will have the
same value on $r^n$ as the map $g$ from 3.4.1. Checking that if
$f(\pm 1_n)=(n_1,...,n_k)=h(\pm 1_l)$ then two definition of $\phi
(n_1,...,n_k)$ agree goes essentially the same way as in the proof
of 2.2 and is left to the reader. Similarly one can show the
multiplicativity of $\phi$.

\bigskip

{\bf Remark 3.5:} Any multiplicative map $\phi: H{\Z}\to R$
determines a formal difference law in $R$. It is given by the
formula $r=\phi (\pm 1_2)$. Hence we can say, that the set of
formal difference  laws in $R$ is in natural bijection with the
set of multiplicative maps $H\Z \to R$.

\bigskip

Perhaps it is now a good point to present how our definition works
in known cases, for example in the case of the spectrum $DB$.  In
Section 2 we mentioned that every  formal sum law in $DB$
determines the formal group law in the ordinary sense. Observe now
that a formal difference law $r\in R[2]$ in the sense of
definition 3.1 determines its sum version $w\in R[2]$ by the
formula
$$w=p^3_2\circ p^4_2 (r^2).$$

\noindent Moreover the map $H\N \to R$ defined by $w$ factors
through the map $H\Z \to R$ defined by $r$.

As another example we would like to show how our definition works
in the case of endomorphism $\g$-ring. This notion is probably
less known so we sketch the definition of it following the
presentation from [S, 13.3].

\medskip

{\bf Example:} Let $\cal C$ be a category with $0$-object and
finite coproducts. The natural enrichment of $\cal C$ over $\ga$
is given by
$$X\wedge [k]\ =\ X  \sqcup ... \sqcup X\ \ \ \  (k-fold\ \
coproduct).$$

\noindent Every object $X$ in $\cal C$ has its endomorphism
$\g$-ring denoted $End_{\cal C}(X)$ defined by

$$End_{\cal C}(X)([k])=Hom_{\cal C}(X,X\wedge [k]).$$

\noindent The unit map ${\bf S}\to End_{\cal C}(X)$ comes from the
identity map in $End_{\cal C}(X)([1])$. The multiplication is
induced by the composition product
$$End_{\cal C}([k])(X)\wedge End_{\cal C}(X)([l]) \to End_{\cal C}(X)([k]\wedge [l])$$

$$f\wedge g \mapsto (f\wedge [l])\circ g.$$

\noindent As Schwede points out, every abelian cogroup object
structure on  $X$ determines the map $H{\Z} \to End_{\cal C}(X)$
defined as follows. At a finite pointed set $[k]$ the map

$$H{\Z} ([k])=\tilde{\Z} [k]\to Hom_{\cal C}(X,X\wedge [k])$$

\noindent is an additive extension of the map which sends $i\in
[k]$ to  the $i$th coproduct inclusion $X\to X\sqcup ...\sqcup X$.

Observe now, that every formal difference  law $r\in End_{\cal
C}(X)[2]=Hom_{\cal C}(X,X\sqcup X)$ defines the abelian cogroup
structure on $X$. The co-addition is given by a sum version of $r$

 $$p^3_2\circ p^4_2 (r^2)\in Hom_{\cal C}(X,X\sqcup X).$$

 \noindent It is abelian because of the invariance of formal
 sum laws under the action by permutations. By the same reason
 the associativity condition is fulfilled. The co-inverse is given
 by $p^2_1(r)\in End_{\cal C}(X)([1])=Hom_{\cal C}(X,X)$. The
 co-unit equals $s^2_{1,2,1}(r)\in End_{\cal C}(X)([1])=Hom_{\cal
 C}(X,X)$. The described by Schwede (and recalled above) map $H{\Z}\to End_{\cal C}(X)$
 agrees with the one obtained by theorem 3.4 from $r$.

\medskip

We suggest the reader to work out by himself, how our theory works
in the case of matrix  $\g$-rings, see [S, 13.5]. Below we come
back to the question when two formal difference  laws define
homotopic maps $H{\Z} \to R$.

\bigskip

{\bf Theorem 3.6}. Two strictly isomorphic formal difference laws
determine homotopic maps in the space of  maps $H{\Z}\to R$.

\bigskip

Recall that when $a\in R[1]$ then multiplication by $a$ from the
left or from the right  determines the map $m_a:R\to R$. Because
left and right multiplications  are formally the same we will
assume that $m_a$ comes from the multiplication from the left. The
theorem 3.6 follows easily from the following lemma.

\bigskip

{\bf Lemma 3.7}. Assume that $a$ and $b$ are two elements in
$R[1]$ which determine the same element in $\pi_0(R)$ under the
obvious map $R[1]\to R$. Then multiplication maps $m_a$ and $m_b$
are homotopic.

\medskip

Proof. The statement of the lemma follows directly from the
definitions, if one carefully investigates  what does it mean that
$a\in R[1]$ determines an element in $\pi_0(R)$. Observe that
choosing $a\in R[1]$ we uniquely choose a map $f_a:{\bf S}\to R$:
it is fully described on the set $[1]$ where we put ${\bf S}[1]\ni
1\mapsto a\in R[1]$. In higher degrees our map is determined by
this data because every $i\in {\bf S}[n]$ can be viewed as the
image of the map $[1]\to [n]$ taking $1$ to $i$.

Moreover observe that on the set $[1]$ our map can be described as
the unit map $\eta$ multiplied from the left by $a$. By naturality
of the multiplication map we see that our map $f_a$ is equal to
$\eta$ composed with the multiplication from the left by $a$.
Observe now that we can decompose the map $m_a$ as a composition

$$R\to {\bf S}\wedge R\to R\wedge R \to R$$

where the first map is an isomorphism, the second is given by
$f_a\wedge id$ and the third is given by the multiplicative
structure $\mu$ of the $\g$-ring $R$.

 Now we can come back to the proof of 3.7. By the assumption $f_a$
 and $f_b$ determine the homotopic maps of spectra. From this we
 get that $f_a\wedge id$ is homotopic to $f_b\wedge id$ and hence
 $m_a$ and $m_b$ give us the homotopic maps of spectra.

 \bigskip

 Let come back to the proof of 3.6. We know that $r_1$ and
 $r_2$ are strict isomorphic, and the isomorphism is given
 by an invertible element $a\in R[1]$. Let $\phi_1$ ($\phi_2$)
 denote the map $H\Z \to R$ determined by $r_1$ ($r_2$). Then
 $$\phi _2=m_{a^{-1}}\circ \phi_1 \circ m_a$$
 \noindent By the assumption
 $m_a$ and $m_{a^{-1}}$ are homotopic to the $m_1$, hence to the
 identity map. This finishes the proof of 3.6.

\bigskip

The referee suggested the following interesting generalization  of
the considerations above to the case when $R$ is not  discrete. In
the latter case, for every natural $n$, the $n$th simplicial
degree of $R[K]$ assemble to a discrete $\g$-ring $R_n$. Hence, in
the case of $R$ not being discrete,  we can talk about simplicial
set $FDL(R)$ of formal difference laws in $R$ which in degree $n$
has the set of formal difference laws in $R_n$. Similarly we can
consider the simplicial set $\g (H{\Z} , R)$ of multiplicative
maps from $H\Z$ to $R$ which in degree $n$ has the set of such
maps to $R_n$. Then the theorem 3.5 can be stated as

\medskip

{\bf Theorem 3.8:} Simplicial sets $FDL(R)_*$ and $\g (H\Z ,R)_*$
are naturally isomorphic.

\medskip

Moreover the story goes further taking into account  the action of
the invertible elements of $R[1]$. Let $G_n$ be the group of
invertible elements in $R_n[1]$.  they assemble to a simplicial
group $G_*$ and this simplicial group acts on both simplicial sets
from 3.8 by conjugation. With all this structure in mind we can
generalize 3.6 to the following statement.

\medskip

{\bf Theorem 3.9:} The homotopy orbit sets $FDL(R)_{hG_*}$ and $\g
(H{\Z} ,R)_{hG_*}$ are isomorphic.

 \bigskip

 Now we would like to comment a little on the homotopical meaning of our constructions.
 As one can see, the proof of 3.6 was derived directly from the definitions. But
 of course we would like to know whether  two strictly isomorphic
 formal difference  laws define homotopic maps in the space of
 multiplicative maps from $H{\Z} \to R$ or, equivalently, the same element in the
 $0$th homotopy group of the space of multiplicative maps $H{\Z} \to R$ as in [S].
 The answer here is not easy to achieve or even to formulate the
conjecture. Our constructions depended heavily on the small model
of $H\Z$ which is not cofibrant. Schwede's homotopical
calculations were  possible also because of the definition of
$DB$-spectrum and its closed relations to symmetric algebra. With
the lack of all these structures we can only propose the following
weak homotopical statement:

\bigskip

{\bf Proposition 3.10:} Let $r_1$ and $r_2$ be two strictly
isomorphic formal difference laws in a $\g$-ring $R$. There exists
a weak equivalence of $\g$-rings $h:R\to R_3$ such that maps
defined by $r_1$ and $r_2$ composed with $h$ are homotopic in the
space of multiplicative maps $H{\Z} \to R_3$

\medskip
Proof. We will be sketchy here because the proof is taken directly
from [S].  For any $\g$-ring $R$ invertible elements in $R[1]$
which maps to the unit component of $R$  form a group $G$ which
acts by conjugation on $R[2]$ and in general on $R$. Two formal
group laws $r_1,\ r_2\in R[2]$ are strictly isomorphic if they are
in the same orbit of this action. Our problem would be solved if
we could extend our conjugation action described above to the
action of the whole unit component of $R$.

Following [S, section 3] we first choose  $R^f$ to be a stably
fibrant replacement of $R$ in a correct model category structure.
Then we define homotopy units $R^*$ as the union of invertible
components of the simplicial monoid $R^f[1]$. We have
$\pi_0R^*=units (\pi_0R)$ and $\pi_iR^*=\pi_iR$ for $i\geq 1$. The
stable equivalence  $R\to R^f$ gives us a homomorphism $\phi:G\to
R^f[1]$ of simplicial monoids with the image in $R^*$. We want to
extend the conjugation action of $G$ to $R^*$. The problem is that
the conjugation action uses strict inverses while $R^*$ is only a
group-like simplicial monoid. The construction how to get around
this difficulty  goes in several steps (see [S, section 4] for the
details).
\medskip

{\bf Step 1}. We factor the map $G\to R^*$ into  $G\to cR^*\to
R^*$ in the correct model category structure of simplicial monoids
where the first map is a cofibration and the second  an acyclic
fibration . Let $UR^*$ denote the group completion of $cR^*$. Then
$UR^*$ is a simplicial group and Lemma 4.3 of [S] tells us that
the map $cR^*\to UR^*$ is a weak equivalence.

{\bf Step 2}. Let ${\bf S}[cR^*]$ be the monoid $\g$-ring with
coefficients in the sphere spectrum. We take the obvious map ${\bf
S}[cR^*]\to R^f$ and factor it in the model category of $\g$-rings
as a cofibration followed by an acyclic fibration
$${\bf S}[cR^*]\to R_1\to R^f.$$ \noindent Then we define another
$\g$-ring $R_2$ as a pushout, in the category of $\g$-rings of

$$\matrix{ & {\bf S}[cR^*]& \to & R_1  \cr
            &\downarrow & \  & \downarrow  \cr
           &{\bf S}[UR^*] &\to  & R_2    \cr}
  $$

\noindent Lemma 4.4 of [S] tells us that the map $R_1\to R_2$ is a
stable equivalence.

{\bf Step 3}. Now  we define $R_3$ to be a stably fibrant
replacement of $R_2$. The induced map ${\bf S}[UR^*] \to R_3$
induces a weak equivalence between  $UR^*$ and the invertible
components of $R_3[1]$. The simplicial group $UR^*$ acts by
conjugation on $R_3$ via homomorphisms of $\g$-rings and this
action extends  the action of $G$.

\bigskip

{\it Final remark}. Of course it is very tempting to speculate
that the weak homotopy type of the  full space of multiplicative
maps $H{\Z} \to R$ should be described via the classifying space
of the groupoid  of formal difference  laws and strict
isomorphism, as it is proved in [S] in the case of the spectrum
$DB$. So far we do not see any way of attacking this problem in
full generality.

\bigskip

{\bf References:}

\medskip

[D] W. G. Dwyer. {\it Homotopy operations for simplicial
commutative algebras}. Trans. Amer. Math. Soc. 260 (1980), no. 2,
421--435.

[L] M. Lydakis. {\it Smash products and $\Gamma$-spaces}. Math.
Proc. Cambridge Philos. Soc. 126 (1999), no. 2, 311--328.

[S] S. Schwede.{\it Formal groups and stable homotopy of
commutative rings}. Geom. Topol. 8 (2004), 335--412.

[S1] S. Schwede. {\it Stable homotopy of algebraic theories}.
Topology 40 (2001) 1-41.

[S2] S. Schwede. {\it Stable homotopical algebra and $\g$-spaces}.
Math. Proc. Camb. Phil. Soc 126 (1999) 329-356.

[Se] G. Segal.{\it Categories and cohomology theories}. Topology
13 (1974) 293--312.

[Sp] E. Spanier. {\it Infinite symmetric products, function
spaces, and duality}. Annals of Math. 69 (1959) 142-198.

\bigskip

\bf Instytut Matematyki, University of Warsaw

ul.Banacha 2, 02-097 Warsaw, Poland

e-mail: betley@mimuw.edu.pl

\end